Yurii VOLKOV, DSc (Phys. & Math.), Prof.
ORCID ID: 0000-0002-2270-3407
e-mail: yuriivolkov38@gmail.com
Volodymyr Vynnychenko Central Ukrainian State University, Kropyvnytskyi, Ukraine

Oleksandr VOLKOV, Master's Student
ORCID ID: 0009-0004-0247-9921
e-mail: oleksandr_volkov@berkeley.edu
University of California, Berkeley, Berkeley, USA

# ON THE CLASS OF EXPONENTIAL STATISTICAL STRUCTURES OF TYPE B

*The article is devoted to the study of exponential statistical structures of type B, which constitute an important subclass of exponential families of probability distributions. This class is characterized by a number of analytical and probabilistic properties that make it a convenient tool for solving both theoretical and applied problems in mathematical statistics. The relevance of this research lies in the need to generalize known classes of distributions and to develop a unified framework for their analysis, which is essential for applications in stochastic modeling, machine learning, and financial mathematics.*

*The paper proposes a formal definition of type B distributions based on the Laplace transform of dominating measures and a system of functional-differential equations describing their structure. Necessary and sufficient conditions for a statistical structure to belong to class B are established, and it is proved that such structures can be represented through a dominating measure with an explicit Laplace transform. The obtained results make it possible to describe a wide range of well-known one-dimensional and multivariate distributions, including the Binomial, Poisson, Normal, Gamma, Polynomial, and Logarithmic distributions, as well as specific cases such as the Borel–Tanner and Random Walk distributions.*

*Particular attention is given to the proof of structural theorems that determine the stability of class B under linear transformations and the addition of independent random vectors. It is shown that if a distribution belongs to class B, its linear transformations and sums also belong to this class with the corresponding parameters. Recursive relations for initial and central moments as well as for semi-invariants are obtained, providing an efficient analytical and computational framework for their evaluation.*

*Furthermore, the tails of type B distributions are investigated using the properties of the Laplace transform. As a result, new exponential inequalities for estimating the probabilities of large deviations are derived, which extend classical approaches to the analysis of statistical distributions. The obtained results can be applied in theoretical studies and in practical problems of stochastic modeling.*

K e y w o r d s: *exponential statistical structures, class B, probability distributions, Laplace transform, stochastic modeling, linear operators.*

**Introduction**

The concept of an exponential statistical structure was formed in the mid-20th century within the framework of the theory of exponential families of distributions. The first systematic constructions appeared in the works of J.-R. Barra (Barra, 1981), who developed the analytical foundations of mathematical statistics and introduced the concept of natural parameters of the exponential family. Classical discrete distributions — Binomial, Poisson, and Negative Binomial — were generalized in the monograph by Johnson and Kotz (Johnson & Kotz, 1969). Further development of the theory was achieved through the works of A. Noack (Noack, 1950), who proposed a new class of discrete distributions with the property of closure with respect to the addition operation, and C. Morris (Morris, 1982), who systematized Natural Exponential Families with Quadratic Variance Functions (NEF-QVF). The classical works of Kendall and Stuart (Kendall & Stuart, 1969) presented the foundations of the statistical theory of moments and cumulants, which became the basis for further analytical developments.

At the same time, starting from the 2000s, exponential structures have been considered not only as an analytical tool but also as part of the broader framework of information geometry. In particular, in the works of S. Amari (Amari, 2016), Ay and Jost (Ay, Jost, & Lê, 2017), as well as Kass and Vos (Kass & Vos, 2019), it was shown that exponential families can be described through the Riemannian structure of the parameter space, where the Fisher–Rao metric is defined as $g_{ij}(\theta) = \mathbb{E}_\theta \left[ \frac{\partial \ln f(X;\theta)}{\partial \theta_i} \frac{\partial \ln f(X;\theta)}{\partial \theta_j} \right]$. This approach made it possible to establish a connection between classical exponential distributions, Bayesian models, and the geometric properties of statistical manifolds.

Furthermore, generalizations of the classical NEF–QVF constructions have emerged. In particular, in the works of Letac and Mora (Letac & Mora, 1990), natural exponential families with cubic variance functions (NEF–CVF) were studied, which are described by the relation $V(m) = am^3 + bm^2 + cm + d$, where $V(m)$ is the variance function, and $m$ is the mean. Further generalizations were developed by Hamza and Hassairi (Hamza & Hassairi, 2011), who provided new characterizations of the Letac–Mora class through a Monge–Ampère-type equation for the cumulant function $\kappa(\theta)$: $\det\left(\frac{\partial^2 \kappa(\theta)}{\partial \theta_i \partial \theta_j}\right) = F\left(\frac{\partial \kappa}{\partial \theta_1}, \ldots, \frac{\partial \kappa}{\partial \theta_n}, \theta\right)$, which connects the geometric properties of the parameter space with the analytical characteristics of the NEF.

Recently, Laplace transform methods have been reconsidered as a powerful tool for analyzing the asymptotic properties of exponential models. Katsevich et al. (Katsevich, 2024) developed a rigorous approach to the multidimensional Laplace expansion: $I(\lambda) = \int_{\mathbb{R}^n} e^{-\lambda \varphi(x)} a(x)\, dx \sim e^{-\lambda \varphi(x_0)} \left(\frac{2\pi}{\lambda}\right)^{n/2} \frac{a(x_0)}{\sqrt{\det \varphi''(x_0)}}, \lambda \to \infty$, where $x_0$ is the point of minimum of $\varphi(x)$, which provides precise applicability conditions in high dimensions.

 

Despite the existing results within the NEF–QVF and NEF–CVF frameworks, the class of exponential structures defined through the Laplace transform of dominating measures that satisfy a functional-differential equation linking parameter variables and moment characteristics has not yet been systematically described. This gap is filled by the present study.

The object of research is exponential statistical structures of type B, which form a subclass of exponential families of distributions.

The aim and objectives of the research are to establish the necessary and sufficient conditions for a statistical structure to belong to class B, to construct the corresponding analytical apparatus through the Laplace transform of dominating measures, and to prove structural theorems regarding the closure of this class with respect to linear transformations and convolutions.

**1. Exponential statistical structures of type B**

Among exponential statistical structures (Noack, 1950, p. 177–190), one can distinguish a subclass characterized by a number of important analytical properties.

*Definition* 1. We say that a statistical structure $(Q(x), x \in X)$ belongs to class $B$ if there exists a function $\varphi(z, x)$, which is the Laplace transform of the measure $Q(x)$, and this function satisfies the equation:

$$\frac{\partial \varphi}{\partial z} + \frac{\partial \varphi}{\partial x} V(x) + x\varphi = 0, \quad \varphi(0, x) = 1 \quad \forall x \in X. \tag{1}$$

The distributions $Q(x)$ themselves will be called distributions of type $B$.

**Theorem 1.** *For a statistical structure $\{Q(x), x \in X\}$ to belong to class B, it is necessary and sufficient that there exists a measure $\mu$ dominating the distributions $Q(x)$ such that*

$$Q(x)(dt) = (u(s(x)))^{-1} \exp(-s(x)t)' \mu(dt),$$

*where $u(z)$ is the Laplace transform of the measure $\mu$, and $s(x)$ is the $s$-characteristic of the measure $\mu$.*

Sufficiency is verified directly. Let us prove necessity.

*Proof.* Since $\varphi(0, x) = 1$, the function $\varphi(z, x)$ is the unique solution of equation (1). Let us express it differently. To do this, fix an arbitrary point $(z, x)$ and denote by $s(x)$ a particular solution of the equation $ds/dx = -W(x) = -(V(x))^{-1}$. Since for all $x \in X$ the matrix $W(x)$ is positive definite, there exists a neighborhood $\delta$ of the point $x$ in which the mapping $s(x)$ has an inverse. Denote it by $x(s)$, $s \in \delta$. Next, consider the equation $db/dx = x(s)$. Since $dx_i/dx_j = dx_j/dx_i, i, j = \overline{1, m}$ (as a consequence of the symmetry of the matrix $W(x)$), this equation has at least one solution. Denote by $b(s)$ the analytic continuation of this solution to the strip $T_\delta = \{\zeta \in \mathbb{C}^m : \operatorname{Re} \zeta \in \delta\}$. Then, by direct substitution into equation (1), it can be verified that

$$\varphi(z, x) = \exp(\langle b(s(x)) - b(s(x) + z)\rangle).$$

Consider the measure

$$\mu(dt) = \exp(s(x)t' - b(s(x)))Q(x)(dt).$$

This measure does not depend on $x$, since its Laplace transform is the function $u(z) = \exp(-b(z))$, $z \in T_\delta$. Let us prove that the measure $\mu$ is the desired one. Indeed, consider the measure

$$q(x)(dt) = (u(s(x)))^{-1} \exp(-s(x)t')\mu(dt).$$

The Laplace transform of the measure $q(x)$ is the function $\exp(\langle b(s(x)) - b(s(x) + z)\rangle)$, i.e., $\varphi(z, x)$.

Due to the one-to-one correspondence between measures and their Laplace transforms, this means that $q(x) = Q(x)$, which completes the proof. □

In practice, the membership of a family of distributions $\{Q(x), x \in X\}$ in the class $B$ is established by checking relation (1).

We symbolically denote that the family $\{Q(x), x \in X\}$ belongs to the class $B$ as:

$$Q(x) \in B(x; V(x)).$$

If $\xi \sim Q(x)$, we also write: $\xi \in B(x; V(x))$.

Many well-known families of distributions belong to the class $B$.

*Example 1.* The Binomial distribution $B_p^n \in B(np; np(1-p))$.

*Example 2.* The Poisson distribution $P_\lambda \in B(\lambda; \lambda)$.

*Example 3.* The Negative Binomial distribution

$$NB_p^n \in B(n(1-p)p^{-1}; np^{-2}(1-p)).$$

*Example 4.* The Normal distribution $\Phi_{\alpha, \sigma^2} \in B(\alpha; \sigma^2)$.

*Example 5.* The Gamma distribution $\Gamma_{\alpha,\lambda} \in B(\lambda/\alpha; \lambda/\alpha^2), (\lambda - \text{fixed})$. In particular, the exponential distribution $\Gamma_{\alpha,1} \in B(\alpha^{-1}; \alpha^{-2})$.

*Example 6.* The Multivariate Normal distribution $\Phi_{\alpha,\sigma^2} \in B(\alpha; \sigma^2)$, where $\alpha = (\alpha_1, \ldots, \alpha_m), \sigma = (\sigma_{ij}), i,j = \overline{1,m}$.

*Example 7.* The Multinomial distribution

$$B^n_{(p_1,\ldots,p_m)} \in B(n(p_1,\ldots,p_m); np_i(\delta_{ij} - p_j)), i,j = \overline{1,m}.$$

*Example 8.* The Negative Multinomial distribution (Johnson & Kotz, 1969, p. 292)

$$NB^n_{(p_1,\ldots,p_m)} \in B(n(p_1,\ldots,p_m); np_i(\delta_{ij} + p_j)), i,j = \overline{1,m}.$$

*Example 9.* Consider the Multivariate Logarithmic distribution $\Lambda_{(\theta_1,\ldots,\theta_m)}$ (Johnson & Kotz, 1969, p. 303), defined as follows: let $k = (k_1,\ldots,k_m) \in \mathbb{N}^m, k_1 + \cdots + k_m \geqslant 1$. Then

$$\xi \in \Lambda_{\{\theta_1,\ldots,\theta_m\}} \iff$$
$$\iff \mathsf{P}\{\xi = (k_1,\ldots,k_m)\} = -(k_1 + \cdots + k_m - 1)(k_1! \ldots k_m!)^{-1} \theta_1^{k_1} \ldots \theta_m^{k_m} \log(1 - \theta_1 - \cdots - \theta_m)^{-1},$$
$$0 < \theta_i < 1, \quad \theta_1 + \cdots + \theta_m < 1.$$

This family of distributions belongs to the class $B$:

$$\Lambda_{(\theta_1,\ldots,\theta_m)} \in B\left(-\frac{(\theta_1,\ldots,\theta_m)}{(1-|\theta|)\log(1-|\theta|)} \; ; \; -\frac{\theta_i}{(1-|\theta|)\log(1-|\theta|)}\left(\delta_{ij} + \frac{\theta_j(1+\log(1-|\theta|))}{(1-|\theta|)\log(1-|\theta|)}\right)\right),$$

where $(i,j = \overline{1,m}, \quad |\theta| = \theta_1 + \cdots + \theta_m.)$

In the univariate case:

$$\Lambda_\theta \in B\left(-\frac{\theta}{(1-\theta)\log(1-\theta)}, -\frac{\theta}{(1-\theta)^2\log(1-\theta)}\left(1 + \frac{\theta}{\log(1-\theta)}\right)\right).$$

*Example 10.* Define the Multivariate Random Walk distribution $B^n_{(p_1,\ldots,p_m)}$.

Let

$$k \in \mathbb{N}^m, \quad |k| = k_1 + \cdots + k_m, \quad 0.5 < p_i < 1, \; i = \overline{1,m}, \quad p = \frac{1}{2}\left(1 + \left(\sum_{i=1}^m (2p_i - 1)^{-1}\right)^{-1}\right).$$

Then

$$\xi \in B^n_{(p_1,\ldots,p_m)} \iff \mathsf{P}\{\xi = k\} = \frac{n(|k|-1)}{|k|!}\left(\frac{(|k|+n)/2}{|k|}\right) p^{|k|+n}(1-p)^{|k|-n} \prod_{i=1}^m \left(\frac{2p-1}{2p_i-1}\right)^{k_i}.$$

We have:

$$B^n_{(p_1,\ldots,p_m)} \in B\Big(n((2p_1-1)^{-1},\ldots,n(2p_m-1)^{-1}); n(2p_i-1)^{-1}\big(\delta_{ij} - (2p_j-1)^{-1}\big)\big(2(2p-1) - (2p-1)^{-1}\big)\Big),$$

where $i,j = \overline{1,m}$.

In the univariate case:

$$B^n_p \in \left(n(2p-1)^{-1}; n(2p-1)^{-1}((2p-1)^{-2} - 1)\right).$$

*Example 11.* Define the Multivariate Borel–Tanner distribution $BT^n_{(\alpha_1,\ldots,\alpha_m)} \in B$.

Let

$$k \in \mathbb{N}^m, \quad 0 < \alpha_j < 1, \quad i = \overline{1,m}, \quad \alpha = 1 - \left(\sum_{i=1}^m (1-\alpha_i)^{-1}\right)^{-1}.$$

Then

$$\xi \in BT^n_{(\alpha_1,\ldots,\alpha_m)} \iff \mathsf{P}\{\xi = k\} = \frac{|k|!}{k!(|k|-n)!} n|k|^{|k|-n-1} \prod_{i=1}^m \left(\frac{1-\alpha}{1-\alpha_j}\right)^{k_i} \exp(-|k|\alpha)\alpha^{|k|-n}.$$

We have:

$$BT^n_{(\alpha_1,\ldots,\alpha_m)} \in B\Big(n\big((1-\alpha_1)^{-1},\ldots,(1-\alpha_m)^{-1}\big); \quad n(1-\alpha_j)^{-1}\big(\delta_{ij} - (1-\alpha_j)^{-1}(2 - ((1-\alpha)^{-1} - \alpha))\big)\Big),$$
$$i,j = \overline{1,m}.$$

In the univariate case:

$$BT^n_\alpha \in B\left(\frac{n}{1-\alpha}; \frac{n\alpha}{(1-\alpha)^3}\right).$$

The class $B$ also includes the Noack distribution family [3] and the so-called natural exponential families with quadratic variance function (NEF-QVF distributions) (Morris, 1982).

**Theorem 2.** *Let $a = (a_{ij}), i, j = \overline{1, m}$ be an arbitrary non-singular matrix, and $b = (b_1, \ldots, b_m)$ an arbitrary vector. Then, if $\xi \in B(x, V(x))$, we have*

$$\eta = a\xi + b \in B(ax + b; aV(x)a').$$

*Proof.* To prove the statement, it is sufficient to verify the equality

$$\frac{\partial \varphi_\eta}{\partial z} + \frac{\partial \varphi_\eta}{\partial y} V_\eta(y) - y\varphi_\eta = 0,$$

with $y = ax + b$, where $\varphi_\eta = \varphi_\eta(z, y)$ is the Laplace transform of the distribution of the random vector $\eta$, and $V_\eta(y) = \text{cov}\,\eta$. The verification of the given equality is carried out in coordinate form.

Note that $M_\eta = y = ax + b$, hence

$$y_j = \sum_{k=1}^{m} a_{jk} x_k + b_j, \quad j = \overline{1, m}, \quad \text{cov}\,\eta = \text{cov}(a\xi + b) = aV(x)a',$$

therefore

$$V_\eta = \left( \sum_{k,r=1}^{m} a_{ir} a_{jr} \nu_{rk}(x) \right), \quad i, j = \overline{1, m}, \quad \varphi_\eta = M \exp(-z(a\xi + b)') = \exp(-zb')\varphi(az, x).$$

Let $\omega = az$, $\omega = (\omega_1, \ldots, \omega_m)$. Then

$$\frac{\partial \varphi_\eta}{\partial z_i} + \sum_{j=1}^{m} \left( \sum_{k,r=1}^{m} a_{ir} a_{jr} \nu_{rk}(x) \right) \frac{\partial \varphi_\eta}{\partial y_j} + y_j \varphi_\eta =$$

$$= \exp(-zb') \left( \sum_{r=1}^{m} a_{ri} \frac{\partial \varphi(\omega, x)}{\partial \omega_r} + \sum_{q=1}^{m} \frac{\partial \varphi(\omega, x)}{\partial x_q} \sum_{k,r=1}^{m} \left( \sum_{j=1}^{m} a_{qj}^{(-1)} a_{jk} \right) a_{ir} \nu_{rk}(x) + \sum_{r=1}^{m} a_{ir} x_r \varphi(\omega, x) \right),$$

where $a_{qj}^{(-1)}$ are the elements of the matrix $a^{-1}$.
Since

$$\sum_{j=1}^{m} a_{qj}^{(-1)} a_{jk} = \delta_{qk},$$

we obtain

$$\exp(-zb') \sum_{r=1}^{m} a_{ir} \left( \frac{\partial \varphi(\omega, x)}{\partial \omega_r} + \sum_{q=1}^{m} \frac{\partial \varphi(\omega, x)}{\partial x_q} \nu_{rq}(x) + x_r \varphi(\omega, x) \right) = 0,$$

which completes the proof. $\square$

**Theorem 3.** *Let $\xi_1, \ldots, \xi_n, \ldots$ be a sequence of independent and identically distributed random vectors, with $\xi_1 \in B(x; V(x))$. Then*

$$\eta_n = \xi_1 + \cdots + \xi_n \in B(nx; nV(x)).$$

The validity of this statement follows from the fact that

$$\varphi_{\eta_n}(z, x) = (\varphi(z, x))^n.$$

**Corollary 1.**
$$\zeta_n := (\xi_1 + \cdots + \xi_n)/n \in B(x; n^{-1}V(x)).$$

**Theorem 4.** *If $Q(x) \in B(x; V(x))$, then $I(x) = V^{-1}(x)$, where $I(x)$ is the Fisher information matrix of the family $\{Q(x), x \in X\}$.*

*Proof.* Let the measure $\mu$ dominate the distributions $Q(x)$. Then

$$(dQ(x)/d\mu)(t) = p(t, x) = \exp(-s(x)t' - \ln u(s(x))).$$

Let $q_{ij}(x), i, j = \overline{1, m}$, denote the elements of the matrix $I(x)$, and let $\xi \sim Q(x)$. Then

$$q_{ij}(x) = M \left( \frac{\partial \ln p(\xi, x)}{\partial x_i} \cdot \frac{\partial \ln p(\xi, x)}{\partial x_{i,j}} \right) = M \left( (p(\xi, x))^{-2} \frac{\partial p(\xi, x)}{\partial x_i} \cdot \frac{\partial p(\xi, x)}{\partial x_{i,j}} \right),$$

and if we denote $\omega := s(x)t' + \ln u(s(x))$, then

$$q_{ij}(x) = L_1 \left( e^{2\omega} (e^{-\omega})_x^{(e_i)} (e^{-\omega})_x^{(e_j)}; x \right) = \omega_{ij}(x) = W(x) = V^{-1}(x),$$

which completes the proof. □

**Theorem 5.** *Let $Q(x) \in B(x; V(x))$, $a_k = a_k(x)$ be the raw moments of order $k = (k_1, \ldots, k_m)$ of the distribution $Q(x)$, $\beta_k = \beta_k(x)$ the central moments, and $\sigma_k = \sigma_k(x)$ the semivariances. Then for all $k \in \mathbb{N}^m$ the following recurrence relations hold:*

$$a_{k+e_i} = \sum_{j=1}^m \nu_{ij}(x) \frac{\partial a_k}{\partial x_j} + x_i a_k, \quad a_{e_i} = x_i, \quad i = \overline{1,m}; \tag{2}$$

$$\beta_{k+e_i} = \sum_{j=1}^m \nu_{ij}(x) \left( \frac{\partial \beta_k}{\partial x_j} + k_j \beta_{k-e_j} \right), \quad \beta_0 = 1, \quad \beta_{e_i} = 0, \quad i = \overline{1,m}; \tag{3}$$

$$\sigma_{k+e_i} = \sum_{j=1}^m \nu_{ij}(x) \frac{\partial \sigma_k}{\partial x_j}, \quad \sigma_{e_i} = x_i, \quad i = \overline{1,m}. \tag{4}$$

*Proof.* Rewrite relation (1) in coordinate form:

$$\frac{\partial \varphi}{\partial z_i} + \sum_{j=1}^m \nu_{ij}(x) \frac{\partial \varphi}{\partial x_j} + x_i \varphi = 0, \quad i = \overline{1,m}. \tag{5}$$

Since $a_k = (-1)^k \varphi^{(k)}(0, x)$, differentiating equality (5) $k$ times and setting $z = 0$ yields (2).

Let $\beta_k(x)$ be the central moments related to the initial moments $a_k(x)$ by the relation

$$\beta_k(x) = \sum_{p \leqslant k} (-1)^{|k-p|} C_k^p \, a_p(x) \, a_1(x)^{k-p}.$$

Differentiating with respect to $x_j$ and taking into account (2), we have

$$\frac{\partial \beta_k}{\partial x_j} = \sum_{p \leqslant k} \frac{\partial \beta_k}{\partial a_p} \left( a_{p+e_j} - x_j a_p \right) = \sum_{p \leqslant k} \frac{\partial \beta_k}{\partial a_p} a_{p+e_j} - k_j \beta_{k-e_j}.$$

Substituting this into (2), we obtain (3).

To obtain relations (4), note that $\sigma_k = (-1)^k (\ln \varphi(z,x))^{(k)}_{z=0}$. Moreover, from relations (1) it follows that the function $\psi = \ln \varphi(z, x)$ satisfies the equation

$$\frac{\partial \psi}{\partial z} + V(x) \frac{\partial \psi}{\partial x} + x = 0,$$

or in coordinate form

$$\frac{\partial \psi}{\partial z_i} + \sum_{j=1}^m \nu_{ij}(x) \frac{\partial \psi}{\partial x_j} + x_i = 0, \quad i = \overline{1,m},$$

from which (4) follows. □

**Corollary 2.** *In the one-dimensional case $(m = 1)$ the following relations hold:*

$$a_{k+1} = V(x) \frac{\partial a_k}{\partial x} + x a_k, \quad a_1 = x, \quad k = 1, 2, \ldots$$

$$\beta_{k+1} = V(x) \left( \frac{\partial \beta_k}{\partial x} + k \beta_{k-1} \right), \quad \beta_0 = 1, \quad \beta_1 = 0, \quad k = 1, 2, \ldots$$

$$\sigma_{k+1} = V(x) \frac{\partial \sigma_k}{\partial x}, \quad \sigma_1 = x, \quad k = 1, 2, \ldots \tag{6}$$

The obtained recurrence relations for moments and semivariances are quite convenient for their computation. In some cases, all moments and semivariances can be derived explicitly from these relations.

We limit ourselves to some examples for finding semivariances, since there exist formulas expressing moments in terms of semivariances (Kendall & Stuart, 1969).

*Example 12.* Let $V(x) = x(ax + b)$, where $a, b$ are arbitrary constants. Then the solution of equation (6) (can be verified directly) is the function

$$\sigma_{k+1}(x) = x(ax + b) \sum_{m=1}^k m! \, b^{k-m} S(k,m)(ax)^{m-1}, \quad k = 1, 2, \ldots \tag{7}$$

where $S(k,m)$ are the Stirling numbers of the second kind, i.e.,

$$S(k,m) = \frac{1}{m!} \sum_{j=0}^m \binom{m}{j} (-1)^j (m-j)^k.$$

In particular, setting $a = -1$, $b = 1$, $x = np$ in (7), we obtain explicit formulas for the semivariances of the binomial distribution $B_p^n$:

$$\sigma_{k+1} = np(1-p) \sum_{m=1}^{k} m!\, S(k,m) p^{m-1}, \quad k = 1, 2, \ldots$$

setting $a = 1$, $b = 1$, $x = n(1-p)p^{-1}$, we obtain explicit formulas for the semivariances of the negative binomial distribution $NB_p^n$:

$$\sigma_{k+1} = np^{-2}(1-p) \sum_{m=1}^{k} m!\, S(k,m)\, p^{1-m}(1-p)^{m-1}, \quad k = 1, 2, \ldots$$

setting $a = 0$, $b = 1$, $x = \lambda$, we obtain the semivariances of the Poisson distribution $P_\lambda$:

$$\sigma_k = \lambda, \quad k = 1, 2, \ldots \tag{8}$$

setting $a = \lambda^{-1}$, $b = 0$, $x = \lambda a^{-1}$, we obtain the semivariances of the gamma distribution $\Gamma_{\lambda, a}$:

$$\sigma_k = (k-1)!\, \lambda a^{-k}. \tag{9}$$

Relations (8) and (9) are known (Kendall & Stuart, 1969, p. 108–109).

*Example 13.* Consider Bernoulli random walks on a line: a particle moves along the axis $t$ at times $t = 1, 2, \ldots$, it moves left with probability $p > 1/2$ and right with probability $q = 1 - p$ by one step. Suppose initially the particle is at position $t = n$. The walk ends if the particle reaches the origin. Let $\xi$ be the random variable representing the time of the walk. This random variable can take values $n, n+2, n+4, \ldots$ with probabilities

$$\omega_{nk} = \frac{n}{k} \binom{k}{(k+n)/2} p^{(k+n)/2} q^{(k-n)/2}, \quad k = n, n+2, \ldots \quad (Volkov, 1992, p.152).$$

Since the generating function of this distribution is $u(s) = \left(1 - \sqrt{1 - 4pqs^2}\right)^n (2qs)^{-n}$, it can be easily verified that

$$\xi \in B\left(\frac{n}{2p-1}, \left(\frac{n}{2p-1}\right)^3 - \frac{n}{2p-1}\right),$$

thus the semivariances of the random variable $\xi$ satisfy the equation:

$$\sigma_{k+1}(x) = (x^3 n^{-2} - x)\sigma_k'(x), \quad \sigma_1(x) = x, \quad x = \frac{n}{2p-1}, \quad k = 1, 2, \ldots$$

It can be directly verified that the solution of this equation is the function

$$\sigma_{k+2} = (x^3 n^{-2} - x) \sum_{r=0}^{m} (-1)^{k+r}(2r+1)!! c_{k,r}(xn^{-1})^{2r},$$

where

$$c_{k,r} = \sum_{m=r}^{k} \binom{k}{m} 2^{m-r} S(m,r), \quad k = 0, 1, 2, \ldots$$

Next, we will consider the estimation of the "tails" of distributions of type $B$.

**Theorem 6.** *Let $Q(x) \in B(x; V(x))$, $x \in X$, and $s(x)$ be the $s$-characteristic of the measure $\mu$ dominating the distributions $Q(x)$. Then for any $y$ such that $s(x) \geqslant s(y)$, the following inequality holds:*

$$Q(x)(\{t : t \geqslant y\}) \leqslant \exp\left(-(y-x)\left(\int_0^1 (1-t)V'(x+t(y-x))dt\right)(y-x)'\right). \tag{10}$$

*Remark* 1. If $t = (t_1, \ldots, t_m)$, $y = (y_1, \ldots, y_m)$, then

$$t \geqslant y \Leftrightarrow t_1 \geqslant y_1, \ldots, t_m \geqslant y_m;$$

$$s(x) \geqslant s(y) \Leftrightarrow s_1(x_1, \ldots, x_m) \geqslant s_1(y_1, \ldots, y_m), \ldots, s_m(x_1, \ldots, x_m) \geqslant s_m(y_1, \ldots, y_m).$$

*Proof.* Consider the operator

$$L_1(f; x) = \int_{\mathbb{R}^m} f(t) Q(x)(dt).$$

Let $\mu$ dominate the distributions $Q(x)$, $s(x)$ be the $s$-characteristic of $\mu$, and $u(z)$ be the Laplace transform of $\mu$. Then, if $\eta(t)$ is the Heaviside unit function, for any vector $a = (a_1, \ldots, a_m)$ with non-negative coordinates, $\eta(t_1) \ldots \eta(t_m) \leqslant \exp(at')$, and hence for any $y \in \mathbb{R}^m$

$$L_1\big(\eta(t_1-y_1)\ldots\eta(t_m-y_m); x\big) \leqslant L_1\big(\exp(a(t-y)'); x\big) = \exp\left(-\left(ay' - \ln\frac{u(s(x)-a)}{u(s(x))}\right)\right).$$

The last inequality is equivalent to

$$Q(x)\{t : t \geqslant y\} \leqslant \exp\left(-\left(ay' - \ln\frac{u(s(x)-a)}{u(s(x))}\right)\right). \tag{11}$$

Set $\beta(a) = ay' - \ln\frac{u(s(x)-a)}{u(s(x))}$, and using properties of the Laplace transform of $\mu$, we obtain

$$\frac{d\beta}{da} = y - x(s(x) - a), \quad x(s) = -d\ln u(s)/ds$$

– the mapping inverse to the mapping

$$s(x): \frac{d^2\beta}{da^2} = -V(x(s(x)-a))).$$

Since the matrix $V$ is positive definite, the function $\beta(a)$ is convex, and hence at the point $a(y) = (a_1(y), \ldots, a_m(y))$ defined by the system

$$y = x(s(x) - a).$$

it attains a maximum. Therefore,

$$a(y) = s(x) - s(y).$$

Set

$$A(y) = \beta(a(y)) = a(y)y' - \ln\frac{u(s(x)-a(y))}{u(s(x))},$$

which is the maximum value of the function $\beta(a)$. The expression for $A(y)$ can be simplified. Note that $A(x) = 0$. Moreover,

$$\frac{dA}{dy} = \left(\frac{da_1}{dy}, \ldots, \frac{da_m}{dy}\right)y' + a(y) - x(s(y))\left(\frac{ds_1}{dy}, \ldots, \frac{ds_m}{dy}\right) = a(y).$$

Hence, $\frac{dA}{dy}\Big|_{y=x=0}$. Further,

$$\frac{d^2A}{dy^2} = \frac{da}{dy} = -\frac{ds(y)}{dy} = -V^{-1}(y).$$

Then, representing the function $A(y)$ by its Taylor expansion around the point $x$ with the remainder in integral form, we obtain

$$A(y) = (y-x)\left(\int_0^1 (1-t)V^{-1}(x+t(y-x))dt\right)(y-x)'.$$

Since inequality (11) holds for any $a$ with non-negative coordinates, it also holds for $a = a(y)$, if $s(x) \geqslant s(y)$ (for example, if all elements of $V^{-1}$ are non-negative, then $a(y) \geqslant 0$ provided $x \leqslant y$). Hence, $Q(x)\{t : t \geqslant y\} \leqslant \exp(-A(y))$, which completes the proof. □

**Discussion and conclusions**

The relations (2)–(4) obtained in this work generalize the classical results of Kendall and Stuart (Kendall & Stuart, 1969) regarding moment equations, but unlike them, they are derived within a unified approach to exponential structures defined through the Laplace transform. This allows for the derivation of recurrent formulas for raw moments, central moments, and semi-invariants not only for individual distributions but also for a wide class of measure-dominated families that satisfy the functional-differential equation (1). Such an approach provides a unified algorithm for computing analytical characteristics for various distributions, including Binomial, Poisson, Normal, Gamma, and Lognormal, which turned out to be special cases of the class $B(x; V(x))$.

Of particular importance are the obtained structural theorems (2)–(4), which prove the closure of the class $B$ with respect to linear transformations and convolutions. This means that when independent random vectors are added or affine transformations of the distribution parameters are applied, the distribution remains within the same class. This property is fundamental for constructing multivariate stochastic models, particularly in problems of risk aggregation and modeling of cumulative random effects in financial processes.

An important result of the research is also the proof of exponential inequalities (10), which provide estimates of large deviation probabilities. This opens up the possibility of applying the methods of class $B$ in problems of risk control, system reliability assessment, and the analysis of the tail behavior of distributions in statistical models.

The practical significance of the obtained results lies in the fact that structures of type $B$ form a theoretical basis for constructing generalized exponential models capable of describing both symmetric and asymmetric stochastic processes with exponential tails. This makes them promising for modeling loss distributions in insurance, the behavior of financial assets, as well as for training machine learning models that take into account nonlinear variance dependencies.

Further research should be directed towards extending the theory of exponential structures of type $B$ to non-stationary processes, analyzing multivariate parametric models with dependent components, and developing efficient numerical algo-rithms for computing moments and semi-invariants for high-dimensional distributions. A separate promising direction is the application of the obtained results to problems of optimal statistical inference, parameter estimation in complex stochastic systems, and the construction of adaptive models in machine learning and financial analytics.

**References**

Amari, S. (2016). *Information geometry and its applications*. Springer.
Ay, N., Jost, J., & Lê, H. V. (2017). *Information geometry*. Springer.
Barra, J.-R. (1981). *Mathematical basis of statistics*. Academic Press, New York.
Hamza, M., & Hassairi, A. (2011). New characterizations of the letac–mora class of real cubic natural exponential families. *arXiv preprint arXiv:1103.1396* , .
Johnson, N. L., & Kotz, S. (1969). *Distributions in statistics: Discrete distributions*. Houghton Mifflin, Boston.
Kass, R. E., & Vos, P. W. (2019). *Geometrical foundations of asymptotic inference*. Springer.
Katsevich, A. (2024). The laplace asymptotic expansion in high dimensions. *arXiv preprint rXiv:2406.12706* , .
Kendall, M., & Stuart, A. (1969). *The advanced theory of statistics. vol. 1: Distribution theory*. Charles Griffin, London.
Letac, G., & Mora, M. (1990). Natural real exponential families with cubic variance functions. *The Annals of Statistic*, 18(1), 1–37. https://doi.org/10.1214/aos/1176347491
Morris, C. N. (1982). Natural exponential families with quadratic variance functions. *nnals of Mathematical Statistics,*, 10(2), 65—80.
Noack, A. (1950). A class of random variables with discrete distributions. *nnals of Mathematical Statistics,*, 21(1), 127-–132.
Volkov, Y. I. (1992). *Positive operators. approximation. probability*. NMC VO, Kyiv. [in Ukrainian].

**Юрій ВОЛКОВ**, д-р фіз.-мат. наук, проф.
ORCID ID: 0000-0002-2270-3407
e-mail: yuriivolkov38@gmail.com
Центральноукраїнський державний університет імені Володимира Винниченка, Кропивницький, Україна

**Олександр ВОЛКОВ**, магістр
ORCID ID: 0009-0004-0247-9921
e-mail: oleksandr_volkov@berkeley.edu
Університет Каліфорнії, Берклі, Берклі, США

## КЛАС ЕКСПОНЕНЦІАЛЬНИХ СТАТИСТИЧНИХ СТРУКТУР ТИПУ В

*Стаття присвячена дослідженню експоненціальних статистичних структур типу В, які становлять важливий підклас експоненціальних сімей розподілів. Цей клас вирізняється низкою аналітичних і ймовірнісних властивостей, що робить його зручним інструментом для теоретичних і прикладних задач математичної статистики. Актуальність теми зумовлена потребою в узагальненні відомих класів розподілів та побудові єдиного апарату для їх аналізу, що має практичне значення у стохастичному моделюванні, машинному навчанні та фінансовій математиці.*

*У роботі запропоновано формальне означення розподілів типу В на основі перетворення Лапласа домінуючих мір та системи функціонально-диференціальних рівнянь, що описують їх структуру. Встановлено необхідні й достатні умови належності статистичної структури до класу В, доведено, що такі структури можуть бути подані через домінуючу міру з явним перетворенням Лапласа. Отримані результати дозволяють описати широкий спектр відомих розподілів, серед яких біноміальний, пуассонівський, нормальний, гамма-розподіл, поліноміальний, логарифмічний, а також специфічні випадки — розподіл Бореля–Таннера та розподіли випадкових блукань.*

*Особливу увагу приділено доведенню структурних теорем, що визначають стійкість класу В відносно лінійних перетворень та операції додавання незалежних випадкових векторів. Показано, що якщо розподіл належить до класу В, то його лінійні перетворення та суми також залишаються в цьому класі. Отримано рекурентні співвідношення для початкових і центральних моментів, а також для семінваріантів, що забезпечує ефективний апарат для їх аналітичного та чисельного обчислення.*

*Крім того, досліджено властивості "хвостів" розподілів типу В за допомогою характеристик перетворення Лапласа. У результаті виведено нові експоненціальні нерівності для оцінювання ймовірностей великих відхилень, які розширюють класичні підходи до аналізу статистичних розподілів. Отримані результати можуть бути застосовані у теоретичних дослідженнях та в задачах прикладного стохастичного моделювання.*

**К л ю ч о в і  с л о в а:** *експоненціальні статистичні структури, клас В, ймовірнісні розподіли, перетворення Лапласа, стохастичне моделювання, лінійні оператори.*